\documentclass[11pt]{article}
\usepackage{amsmath,amsthm,amssymb}
\usepackage{enumerate}
\usepackage{graphicx}
\usepackage{titlesec}
\usepackage{empheq}
\usepackage{color}
\usepackage{xfrac,bigints}

\newtheorem{Th}{Theorem}[section]
\newtheorem{Lem}[Th]{Lemma} 
\newtheorem{Prop}[Th]{Proposition} 
\newtheorem{Cor}[Th]{Corollary}

\newtheorem{Ex}[Th]{Example}

\titleformat*{\section}{\normalsize\bfseries}
\titleformat*{\subsection}{\normalsize\bfseries}

\newcommand{\argmin}{\mathop{\rm arg~min}\limits}

\def\R{{\mathbb R}}

\def\DS{\displaystyle}

\def\setm{{\hspace{-0.3pt}\setminus \hspace{-0.3pt}}}
\def\sgn{\mbox{\rm sgn}}

\makeatletter 
\long\def\@makefntext#1{\parindent 1em\noindent 
\@hangfrom{\hbox to 1.8em{\hss$^{\@thefnmark}$}}#1}
\makeatother

\begin{document}

\begin{flushright}
\footnotesize{\today}
\end{flushright}
\begin{center}
{\large\bf
  Two-step minimization approach to Sobolev-type inequality\\
  with bounded potential in 1D
}
\end{center}

\begin{center}
  Vina Apriliani$^{\mbox{\scriptsize a,b}}$,
  Masato Kimura$^{\mbox{\scriptsize c}}$,
  Hiroshi Ohtsuka$^{\mbox{\scriptsize c}}$\\~\\
\begin{tabular}{l}
  {\scriptsize a)} Graduate School of Natural Science and Technology, Kanazawa University\\
 {\scriptsize b)} Ar-Raniry State Islamic University Banda Aceh\\
{\scriptsize c)} Faculty of Mathematics and Physics, Kanazawa University
\end{tabular}
\end{center}

\renewcommand{\thefootnote}{\fnsymbol{footnote}}
\footnote[0]{E-mail: vina.apriliani@ar-raniry.ac.id (VA), mkimura@se.kanazawa-u.ac.jp (MK), ohtsuka@se.kanazawa-u.ac.jp (HO)}
\renewcommand{\thefootnote}{\arabic{footnote}}

\vspace{0.5cm}
\begin{abstract}
We present a new method to determine the best constant of the
Sobolev-type embedding in one dimension with a norm including 
a bounded inhomogeneous potential term. This problem is closely connected
to the Green function of the Schr\"odinger operator with inhomogeneous potential.
A minimization problem of a Rayleigh-type quotient
in a Sobolev space gives the best constant of the Sobolev embedding.
We decompose the minimization problem into two sub-minimization
problems and show that the Green function provides the minimizer
of the first minimization problem. 
Then, it enables us to derive a new precise estimate
of the best constant and function for inhomogeneous bounded potential cases.
As applications, we give some examples of the inhomogeneous potential
whose best constant and function of the Sobolev-type embedding are explicitly determined. 
\end{abstract}


\section{Introduction}
\setcounter{equation}{0}
In this paper, $L^p(\R)$ and $W^{1,p}(\R)$
denote the real-valued Lebesgue and Sobolev spaces, and we set $H^1(\R)=W^{1,2}(\R)$.
For every $1\leq p\leq \infty$ and every function $u\in W^{1,p}(\R)$,
it is well known to hold the following Sobolev-type inequality:
\begin{equation}\label{So-ineq}
\|u\|_{\infty}\le C\| u\|_{W^{1,p}(\R)},
\end{equation}
where $\|u\|_\infty =\|u\|_{L^\infty(\R)}$, and
$C=C(p)>0$ is a constant independent of $u$;
see \cite[Theorem 8.8]{Bre11} for example.
Indeed,
we can choose $C(p)=p^{1/p}\leq e^{1/e}$ from an elementary calculation.
We consider a generalization for the case of $p=2$ and
discuss the corresponding best constant.
The Sobolev-type inequality of the form
$\|u\|_{L^q(\Omega)}\le C\| u\|_{W^{1,p}(\Omega)}$ has
been intensively studied even in higher dimensions with $\Omega\subset\R^n$
under various settings in connection with applications in PDE problems,
for example, \cite{CTZ19, CM16, Kam08, WKNTY08}.
In particular, the best constant 
for the case with an inhomogeneous potential is closely connected
with the corresponding Schr\"odinger operator.
However, there has not been much precise analysis of
the best constant of the Sobolev-type inequality when
there is a general inhomogeneous potential term.

Let $V\in L^\infty(\R)$ with $\inf_{x\in\R} V(x) >0$
for simplicity. We introduce a norm on $H^1(\R)=W^{1,2}(\R)$ such as
\[
\|u\|_V^2
:=\int_{-\infty}^\infty (|u'(x)|^2+V(x)|u(x)|^2)\,dx,
\]
which gives an equivalent norm to the usual $\|\cdot\|_{W^{1,2}(\R)}$
from the assumptions on $V(x)$.
We define a Rayleigh-type quotient
\begin{equation}\label{rq}
R(u; V):=\frac{\| u\|_V^2}{\|u\|_{\infty}^2},
\end{equation}
and consider a minimization problem of $R(u;V)$ for $u\in H^1(\R)$.
More precisely, we define
\begin{align}
&m(V):=\inf_{u\in H^1(\R),u\not{\equiv}0} R(u; V),\label{mV}\\
&M(V):=\left\{u\in H^1(\R);~m(V)=
R(u; V),~u\not{\hspace{-0.5ex}\equiv}0
\right\}.\label{MV}
\end{align}
We note that $m(V)^{-1/2}$ gives the best constant of the following
Sobolev-type inequality
\begin{equation}
\label{So-ineq-inhom}
\|u\|_\infty\leq C\|u\|_V,
\end{equation}
which is similar to \eqref{So-ineq}.
In this paper, we call $m(V)^{-1/2}$ the Sobolev best constant of
the inequality \eqref{So-ineq-inhom}.
 
Contrary to the precise and detailed analysis of the Sobolev best constant
and function for the homogeneous case $V\equiv const.$, it is not easy to get the precise 
value of the Sobolev best constant and profile of the minimizer, which achieves
the best constant for inhomogeneous potentials.
However, the significance of inhomogeneous potential cases
is increasing in connection with various applications.
In this research, we focus on one-dimensional and 
inhomogeneous potential cases. This paper aims to develop
a new two-step minimization approach and establish
new criteria to determine the precise value of the Sobolev
best constant and function for \eqref{So-ineq-inhom}.
We clarify the influence of the potential
$V$ to $m(V)$ and the elements of $M(V)$,
and construct concrete examples using the obtained new criteria.

The main idea that we employed in our analysis is the reduction of 
the minimization problem in the whole $H^1(\R)$ to two-step minimizations 
in $K_a$ and in $a\in\R$, where $K_a:=\{u\in H^1(\R);~u(a)=\|u\|_\infty =1\}$.
In the first minimization in $K_a$, the unique minimizer $u_a$ 
coincides with the Green function for the Schr\"odinger operator
$-\frac{d^2}{dx^2}+V(x)$ with the potential function $V(x)$.
This fact closely links the Sobolev best constant problem to the estimate of the Green function
and gives one of our strong motivations. In the second minimization
step among the parameter $a$, we give a necessary and sufficient condition
for the global minimizer regarding the
fundamental solutions of the linear ODE $-\frac{d^2}{dx^2}u+V(x)u=0$.

The outline of the paper is as follows.
In Section~\ref{subsec:2.1}, we
decompose the minimization problem \eqref{mV} to the two sub-minimization problems
and prove the equivalency between the original minimization
problem and our two-step minimization problem.
In Section~\ref{subsec:2.2},
we first establish precise properties
of the minimizer of the first minimization problem in $K_a$
and a relation to the Green function.
Moreover, by extending the right and left tails of the minimizer,
we construct the fundamental solutions of $-\frac{d^2}{dx^2}u+V(x)u=0$
and give upper and lower estimates for them.
Section~\ref{subsec:2.3} is devoted to the second minimization step.
In Theorem~\ref{FFF}, we establish some useful variational formulae for the minimum value
obtained in the first minimization step.
Section~\ref{subsec:2.4} establishes our main theorem, which
gives three new necessary and sufficient conditions 
for the local minimality in the second minimization step.
As an application,
we construct a nontrivial example of the bounded inhomogeneous potential 
for which the Sobolev best constant and function are explicitly specified.

\section{Two-step minimization approach}\label{subsec:2.1}
\setcounter{equation}{0}

Throughout of this paper, we suppose 
\begin{align}\label{bddV}
V\in L^\infty(\R),\quad
0<v_0\le V(x)\le v_1~\mbox{a.e.}~x\in\R,
\end{align}
for some constants $v_0$ and $v_1$.
Then, for $u,~v\in H^1(\R)$, we define
\begin{align*}
  &(u,v)_V:=\int_\R \left(u'(x)v'(x)+V(x)u(x)v(x)\right)\,dx,\quad \| u\|_V:=(u,u)_V^{\frac{1}{2}},
  \quad I(u;V):=\|u\|_V^2,
\end{align*}
where we abbreviate $\frac{du}{dx}$ as $u'$.
We remark that $(u,v)_V$ defines an inner product on $H^1(\R)$, and
the corresponding norm $\|u\|_V$ is equivalent to the norm
of $H^1(\R)$.

In this paper, without loss of generality, 
we suppose that $u\in H^1(\R)$ (or more generally, $u\in W^{1,p}_{loc}(\R)$
for $p\in [1,\infty]$) always satisfies $u\in C^0(\R)$, since
an element of the function space $W^{1,p}_{loc}(\R)$ has a continuous representation
(Theorem~8.8 of \cite{Bre11}). We also remark that
$u\in H^1(\R)$ satisfies $u\in L^\infty(\R)$ and
$\lim_{|x|\to\infty}u(x)=0$ (Theorem~8.8 and Corollary~8.9 of \cite{Bre11}).

For $a\in \R$, we set
\[
K_a:=\{ u\in H^1(\R);~u(a)=\|u\|_\infty =1\},
\]
and define
\[
F(a):=\inf_{u\in K_a} I(u;V).
\]

The following proposition establishes a decomposition principle
for the minimization problem \eqref{mV}.
It gives a foundation for the two-step minimization method
and enables us to capture the detail of the inhomogeneous potential $V$.
\begin{Prop}[decomposition principle]\label{decomposition}
Let $V$ satisfy \eqref{bddV} and set $m(V)$ as \eqref{mV}. Then, we have
\begin{align}\label{infF}
m(V)
=\inf_{a\in\R}F(a).
\end{align}
\end{Prop}
\proof
We define $\tilde{m}(V):=\inf_{a\in\R}F(a)$.
  Let $\{u_n\}_{n=1}^\infty$ be a minimizing sequence to the infimum of \eqref{mV}.
  Without loss of generality, we can assume $u_n\in K_{a_n}$ for some $a_n\in\R$
  and
$m(V)=\lim_{n\to\infty}I(u_n;V)$.
Taking the limit as $n\to\infty$ in $\tilde{m}(V)\le F(a_n)\le I(u_n;V)$,
we obtain $\tilde{m}(V)\le m(V)$.

On the other hand, we choose $\{a_n\}_{n=1}^{\infty}\subset \R$
and $\{u_{n,k}\}_{k=1}^\infty\subset K_{a_n}$ such that 
\[
\tilde{m}(V)=\lim_{n\to\infty}F(a_n),\quad
F(a_n)=\lim_{k\to\infty}I(u_{n,k};V).
\]
Since $m(V)\le I(u_{n,k};V)$, taking the limit as $k\to\infty$, we have
$m(V)\le F(a_n)$. Then, taking the limit as $n\to\infty$, we obtain $m(V)\le \tilde{m}(V)$.
\qed

\section{First minimization step}\label{subsec:2.2}
\setcounter{equation}{0}
We consider the first minimization step:
\begin{align}\label{Fa}
F(a):= \inf_{u\in K_a}\| u\|_V^2.
\end{align}

\begin{Th}\label{th1}
We suppose the condition \eqref{bddV}
and fix $a\in\R$. There exists a unique minimizer $u_a\in K_a$ to \eqref{Fa}, that is,
\begin{align}
  &u_a=\argmin\limits_{u\in K_a}\|u\|_V^2,\label{uaua}\\
  &F(a)= \min_{u\in K_a}\| u\|_V^2=\| u_a\|_V^2\notag ,
\end{align}
and it satisfies the following properties:
\begin{align}
&u_a\in W^{2,\infty}(\R\setm\{a\})\quad \mbox{and}\quad
u_a''(x)=V(x)u_a(x)~(\mbox{a.e.}~ x\in\R\setm \{a\}).\label{p1}\\
&  
e^{-\sqrt{v_1}|x-a|}\le u_a(x) \le e^{-\sqrt{v_0}|x-a|}\quad(x\in\R),\label{p2}\\
&
\frac{v_0}{\sqrt{v_1}}e^{-\sqrt{v_1}|x-a|}\le \sgn (a-x)  u_a'(x)
\le \frac{v_1}{\sqrt{v_0}}e^{-\sqrt{v_0}|x-a|}\quad(x\in\R\setm \{a\}).\label{p3}
\end{align}
In particular, when $V$ is a positive constant (i.e., $v_0=V=v_1$),
then $u_a(x) =e^{-\sqrt{V}|x-a|}$ holds.
\end{Th}
\proof
We define $L_a:=\{u\in H^1(\R);~u(a)=1\}$.
Since $L_a$ is a closed affine subspace of $H^1(\R)$,
there exists a unique $u_a\in L_a$ defined by
$u_a:=\argmin_{u\in L_a}\|u\|_V^2$
and it
is equivalent to the orthogonality condition
(see Theorem~5.2 and Corollary~5.4 of \cite{Bre11}, for example):
\begin{align}\label{ac}
  u_a\in L_a\quad \mbox{and}\quad (u_a,v)_V=0~\mbox{for all}~ v\in H^1(\R)~\mbox{with}~v(a)=0.
\end{align}
It implies that $u_a''=Vu_a$ in ${\cal D}'(\R\setm \{ a\})$ holds.
As \eqref{p2} is proven for $u_a$ below, we will find that
$u_a$ belongs to $K_a$ and satisfies \eqref{uaua}.

Suppose that $V$ is a positive constant function.
Then, a solution of the linear ODE $u''(x)=Vu(x)$ on $\R$
is written in the form:
\begin{align}
u(x)=c_1e^{\sqrt{V}x}+c_2e^{-\sqrt{V}x}\quad (c_1,c_2\in\R).\label{uAB}
\end{align}
Hence, in this case,
the formula $u_a(x) =e^{-\sqrt{V}|x-a|}$ is immediately derived from \eqref{uAB}.

Next, we prove \eqref{p2} for a general potential $V$ with \eqref{bddV}.
We define $\tilde{u}_a(x):=e^{-\sqrt{v_1}|x-a|}$ and 
\begin{align}
  &w(x):=\max \left(u_a(x),\tilde{u}_a(x)\right)\qquad (x\in\R),\label{wx}\\
&\omega :=\{x\in\R;~u_a(x)<\tilde{u}_a(x)\} :\text{open}.\label{omega}
\end{align}
From Theorem~A.1 in Chapter~II of \cite{Kin80}, $w$ satisfies $w\in L_a$ and 
\begin{align*}
  &w(x)=\tilde{u}_a(x),\quad w'(x)=\tilde{u}_a'(x)\quad \mbox{a.e.}~x\in\omega,\\
  &w(x)=u_a(x),\quad w'(x)=u_a'(x)\quad \mbox{a.e.}~x\in\R\setminus\omega.
\end{align*}

We set $v:=u_a-w$, then $v\in H^1(\R)$, $v(a)=0$, and $v\le 0$ holds.
\begin{align*}
\|u_a\|^2_V-\|w\|^2_V
&=(u_a,v)_V+(w,v)_V
=(w,v)_V\\
&=\int_\omega\left(w'(x)v'(x)+V(x)w(x)v(x)\right)\,dx\\
&=\int_\omega\left(\tilde{u}_a'(x)v'(x)+V(x)\tilde{u}_a(x)v(x)\right)\,dx.
\end{align*}
Since $\omega$ is an open set of $\R$ and $a\notin \omega$, it is decomposed to finite or countable
connected components:
$\omega =\cup_{\lambda\in\Lambda}\omega_\lambda$ 
where $\omega_\lambda=(p_\lambda,q_\lambda)$ is an open interval,
and $\Lambda$ is a finite or countable set.
We remark that
$v(p_\lambda)=0$ holds for $p_\lambda\neq -\infty$, 
and $v(q_\lambda)=0$ for $q_\lambda\neq \infty$.
\[
\int_{\omega_\lambda}\tilde{u}_a'(x)v'(x)\,dx
=\big[\tilde{u}_a'(x)v(x)\big]_{p_\lambda}^{q_\lambda}
-\int_{\omega_\lambda}\tilde{u}_a''(x)v(x)\,dx
=-\int_{\omega_\lambda}v_1\tilde{u}_a(x)v(x)\,dx.
\]
This equality is true even when $p_\lambda =-\infty$
or $q_\lambda=\infty$.
So, it implies
\[
\int_\omega \tilde{u}_a'(x)v'(x)\,dx
=\sum_{\lambda\in \Lambda}\int_{\omega_\lambda}\tilde{u}_a'(x)v'(x)\,dx
=-\sum_{\lambda\in \Lambda}\int_{\omega_\lambda}v_1\tilde{u}_a(x)v(x)\,dx
=-\int_\omega v_1\tilde{u}_a(x)v(x)\,dx.
\]
Hence, we obtain
\begin{align*}
\|u_a\|^2_V-\|w\|^2_V
=\int_\omega (V(x)-v_1)\tilde{u}_a(x)v(x)\,dx \ge 0.
\end{align*}
Since $u_a$ is the unique minimizer of \eqref{uaua} and $w\in L_a$,
$w(x)=u_a(x)$ holds for $x\in\R$.
In other words, $\tilde{u}_a(x)\le u_a(x)$ holds for $x\in\R$.
Another inequality of \eqref{p2} is also shown similarly.
From \eqref{p2}, we conclude that $u_a\in K_a$ and \eqref{uaua} holds.

We next prove \eqref{p3}.
Let us consider the left interval $(-\infty,a)$.
From $u_a''=Vu_a\in L^\infty(-\infty,a)$, $u_a'\in C^0(-\infty,a)$
holds.
Since $u_a'\in L^2(-\infty,a)$, there exists
a sequence $\{r_n\}_{n=1}^{\infty}\subset (-\infty, a)$ such that
$r_n\le -n$ and $|u_a'(r_n)|\le \frac{1}{n}$.
Then, for $x<a$, we have 
\begin{align*}
  u_a'(x)=u_a'(r_n)+\int_{r_n}^x u_a''(t)\,dt
  = u_a'(r_n)+\int_{r_n}^x V(t)u_a(t)\,dt.
\end{align*}
Taking the limit as $n\to \infty$, 
\begin{align*}
  u_a'(x)=\int_{-\infty}^x V(t)u_a(t)\,dt\quad \mbox{for}~x\in (-\infty,a)
\end{align*}
holds and, using \eqref{p2}, we obtain
\begin{align*}
&  u_a'(x)\ge v_0\int_{-\infty}^x e^{-\sqrt{v_1}|t-a|}\,dt
  =\frac{v_0}{\sqrt{v_1}}e^{-\sqrt{v_1}|x-a|},\\
&  u_a'(x)\le v_1\int_{-\infty}^x e^{-\sqrt{v_0}|t-a|}\,dt
  =\frac{v_1}{\sqrt{v_0}}e^{-\sqrt{v_0}|x-a|}.
\end{align*}
Similarly, for $x\in (a,\infty)$, we can derive
the inequalities:
\begin{align*}
  -\frac{v_1}{\sqrt{v_0}}e^{-\sqrt{v_0}|x-a|}
  \le u_a'(x)\le -\frac{v_0}{\sqrt{v_1}}e^{-\sqrt{v_1}|x-a|}.
\end{align*}
Finally, $u_a\in W^{2,\infty}(\R\setm\{a\})$ follows from
$u_a''=Vu_a$ in $(-\infty,a)\cup (a,\infty)$ and \eqref{p3}.
\qed
~\\

From Theorem~\ref{th1}, $u_a'(x)$ exists for $x\neq a$, and $u_a'(a)$ does not exist.
However, the left and right derivatives at $x=a$ exist. 
We denote the right and left derivatives of $u_a(x)$ at $x=a$ by
\[
u_a'(a_\pm):=\lim_{h\to \pm 0}\frac{u_a(a+h)-u(a)}{h}.
\]
\begin{Prop}
The Green function for the differential operator $L:=-\frac{d^2}{dx^2}+V$
is given by
\begin{align}\label{Gxy}
G(x,y)=\frac{u_y(x)}{u_y'(y_-)-u_y'(y_+)}=\frac{u_x(y)}{u_x'(x_-)-u_x'(x_+)} \quad (x,y\in\R).
\end{align}
\end{Prop}
\proof
For $x,y\in\R$, we define
\[
g_y(x):=\frac{u_y(x)}{u_y'(y_-)-u_y'(y_+)}.
\]
Then, $g_y\in H^1(\R)$ holds. For all $v\in H^1(\R)$, we have
\begin{align*}
(g_y,v)_V
&=\int_{-\infty}^y\left(g_y'(x)v'(x)+V(x)g_y(x)v(x)\right)dx
+\int_y^{\infty}\left(g_y'(x)v'(x)+V(x)g_y(x)v(x)\right)dx\\
&=\big[ g_y'v \big]_{-\infty}^y +\int_{-\infty}^y \left(-g_y''(x)+V(x)g_y(x)\right) v(x)dx\\
&\quad +\big[ g_y'v \big]_y^{\infty} +\int_y^{\infty} \left(-g_y''(x)+V(x)g_y(x)\right) v(x)dx\\
&=g_y'(y_-)v(y)-g_y'(y_+)v(y)
=v(y)
=\mbox{}_{H^{-1}(\R)}\left\langle
\delta_y,v
\right\rangle_{H^1(\R)},
\end{align*}
where $\delta_y\in H^{-1}(\R)$ denotes the Dirac delta distribution
at $y$.
This equality implies that $Lg_y=\delta_y$ in $H^{-1}(\R)$.
Furthermore, we obtain the second equality of \eqref{Gxy}
from the following equality for $y,z\in\R$:
\begin{align*}
g_z(y)
=\mbox{}_{H^{-1}(\R)}\left\langle
\delta_y,g_z
\right\rangle_{H^1(\R)}
=(g_y,g_z)_V
=(g_z,g_y)_V
=\mbox{}_{H^{-1}(\R)}\left\langle
\delta_z,g_y
\right\rangle_{H^1(\R)}
=g_y(z).
\end{align*}
\qed

\begin{Prop}[comparison principle of $u_a$]
  For $V,\tilde{V}\in L^\infty(\R)$ with $0<v_0\le V(x)\le\tilde{V}(x)$,
  we set $u_a:=\argmin_{u\in K_a}\|u\|_V^2$
  and $\tilde{u}_a:=\argmin_{u\in K_a}\|u\|_{\tilde{V}}^2$,
  then $u_a(x)\ge\tilde{u}_a(x)$ holds for $a\in\R$ and $x\in\R$.
\end{Prop}
\proof
We define $w(x)$ and $\omega$ by \eqref{wx} and \eqref{omega}.
Then, in a similar way to the proof of Theorem~\ref{th1},
we can derive
\begin{align*}
\|u_a\|^2_V-\|w\|^2_V
=\int_\omega (V(x)-\tilde{V}(x))\tilde{u}_a(x)(u_a(x)-w(x))\,dx \ge 0.
\end{align*}
Since $u_a$ is the unique minimizer of \eqref{uaua} and $w\in K_a$,
$w(x)=u_a(x)$ holds for $x\in\R$.
In other words, $\tilde{u}_a(x)\le u_a(x)$ holds for $x\in\R$.
\qed

\begin{Lem}\label{lem-uab}
Suppose \eqref{bddV}.
For $a\le b$, $u_a$ and $u_b$ defined by \eqref{uaua} satisfy
\begin{align*}
  u_a(x)=\frac{u_b(x)}{u_b(a)}\quad (x\le a),\quad
  u_b(x)=\frac{u_a(x)}{u_a(b)}\quad (x\ge a).
\end{align*}
\end{Lem}
\proof
We define
\[
\tilde{u}_a(x):=\left\{
\begin{array}{ll}
  {\DS \frac{u_b(x)}{u_b(a)}}&(x\le a),\\[10pt]
  u_a(x)&(x>a).
\end{array}\right.
\]
Then, $\tilde{u}_a\in K_a$ and, for $v\in H^1(\R)$ with $v(a)=0$,
it satisfies
\begin{align*}
(\tilde{u}_a,v)_V&=
\frac{1}{u_b(a)}\int_{-\infty}^a\left(u_b'(x)v'(x)+V(x)u_b(x)v(x)\right)\,dx\\
&\hspace{20pt}+\int_{-\infty}^a\left(u_a'(x)v'(x)+V(x)u_a(x)v(x)\right)\,dx\\
&=\frac{1}{u_b(a)}\int_{-\infty}^a\left(-u_b''(x)+V(x)u_b(x)\right)v(x)\,dx\\
&\hspace{20pt}+\int_{-\infty}^a\left(-u_a''(x)+V(x)u_a(x)\right)v(x)\,dx=0.
\end{align*}
Since $\tilde{u}_a$ satisfies the orthogonality condition \eqref{ac},
$\tilde{u}_a=u_a$ holds.
In other words, $u_a(x)=u_b(x)/u_b(a)$ holds for $x\le a$.
The other relation, $u_b(x)=u_a(x)/u_a(b)$ for $x\ge a$, is similarly proven.
\qed

From Lemma~\ref{lem-uab}, we have
\begin{align*}
 u_a'(a_-)=\frac{u_b'(a)}{u_b(a)}\quad (a<b),\quad\quad
u_a'(a_+)=\frac{u_b'(a)}{u_b(a)}\quad (b<a).
\end{align*}
These expressions of $u_a'(a_\pm)$ also imply that
\begin{align}\label{ua'}
  [a\to u_a'(a_\pm)]\in W^{1,\infty}_{loc}(\R).
\end{align}

\begin{Th}\label{phipm}
Under the condition \eqref{bddV},
$\varphi_\pm\in W^{2,\infty}_{loc}(\R)$ uniquely exists such that
\begin{align}\label{ivp}
  \begin{cases}
    \varphi''_\pm(x)=V(x)\varphi_\pm(x)\qquad(\text{a.e. }x\in\R),\\
    \varphi_\pm (0)=1,\\
    \DS{\lim_{x\to\pm \infty}\varphi_\pm(x)=0.}
  \end{cases}
\end{align}
Furthermore, they satisfy the following inequalities:
\begin{align}
  &\sqrt{\frac{v_0}{v_1}}e^{-\max (\sqrt{v_0}x,\sqrt{v_1}x)}\le \varphi_+(x)\le
  \sqrt{\frac{v_1}{v_0}}e^{-\min (\sqrt{v_0}x,\sqrt{v_1}x)}\quad (x\in\R),\label{phi+}\\
  &\sqrt{\frac{v_0}{v_1}}e^{\min (\sqrt{v_0}x,\sqrt{v_1}x)}\le \varphi_-(x)\le
  \sqrt{\frac{v_1}{v_0}}e^{\max (\sqrt{v_0}x,\sqrt{v_1}x)}\quad (x\in\R).\label{phi-}
\end{align}
Moreover, 
\begin{align}\label{uax}
u_a(x)=\begin{cases}
\DS{\frac{\varphi_-(x)}{\varphi_-(a)}\qquad (x<a)}\\[10pt]
\DS{\frac{\varphi_+(x)}{\varphi_+(a)}\qquad (x\ge a)}
\end{cases}
\end{align}
holds for $a\in\R$.
\end{Th}
\proof
We set $m_\pm:=u_0'(0_\pm)$. 
We consider the following initial value problem of linear
ordinary differential equation:
\begin{equation}\label{odeV}
  \begin{cases}
    \varphi''(x)=V(x)\varphi (x)\quad \mbox{a.e.}~x\in\R,\\
    \varphi (0)=1,~\varphi'(0)=m_\pm .
  \end{cases}
\end{equation}
Since the coefficient $V$ may not be continuous,
we consider a mild solution for \eqref{odeV}:
\begin{align}\label{mild}
\varphi\in C^1(\R),\quad \varphi(0)=1,\quad
\varphi'(x)=m_\pm +\int_0^x V(t)\varphi (t)\,dt~~(x\in\R).
\end{align}
It is well known that there exists a unique mild solution $\varphi_\pm\in C^1(\R)$
of \eqref{mild} as it is a linear ODE with bounded coefficient $V\in L^\infty(\R)$.
Moreover, from \eqref{mild}, $\varphi_\pm\in  W^{2,\infty}_{loc}(\R)$ also holds.

We can show that the function defined by the right-hand side of \eqref{uax} satisfies the
condition \eqref{ac} for all $a\in\R$, similar to the proof of Lemma~\ref{lem-uab}.
So, we conclude \eqref{uax}. Then, from the estimate \eqref{p2}, we find that
$\varphi_\pm$ satisfy the initial value problems of \eqref{ivp}.

Next, we prove the estimate \eqref{phi-}. Since $\varphi_-(x)=u_0(x)$ for $x\le 0$,
\eqref{phi-} for $x\le 0$ follows from \eqref{p2}.
For $x>0$, we first remark that $\varphi_-(x)>0$ holds from \eqref{mild}.
If $\varphi_-(x_0)=0$ and $\varphi_-(x)>0$ for $0<x<x_0$,
then there should be an $x_1\in (0,x_0)$ such that $\varphi_-'(x_1)<0$,
but it is impossible from \eqref{mild}.
We define $y(x):=\varphi_-'(x)+\sqrt{v_0}\varphi_-(x)$ for $x\ge 0$.
Then, we have
\begin{align*}
  y'(x)=\varphi_-''(x)+\sqrt{v_0}\varphi_-'(x)=\sqrt{v_0}\varphi_-'(x)+V(x)\varphi_-(x)
  \ge   \sqrt{v_0}\varphi_-'(x)+v_0\varphi_-(x)
  =\sqrt{v_0}y(x).
\end{align*}
Solving this differential inequality, we obtain
$y(x)\ge y(0)e^{\sqrt{v_0}x}$. Since $y(0)=m_-+\sqrt{v_0}$,
this is equivalent to
\begin{align*}
  \varphi_-'(x)+\sqrt{v_0}\varphi_-(x)\ge (m_-+\sqrt{v_0})e^{\sqrt{v_0}x}.
\end{align*}
We solve this differential inequality as follows:
\begin{align*}
  \left( e^{\sqrt{v_0}x}\varphi_-(x)\right)'
  =e^{\sqrt{v_0}x}\left( \varphi_-'(x)+\sqrt{v_0}\varphi_-(x)\right)
  \ge (m_-+\sqrt{v_0})e^{2\sqrt{v_0}x}
  =\left( \frac{m_-+\sqrt{v_0}}{2\sqrt{v_0}}e^{2\sqrt{v_0}x}\right)'.
\end{align*}
Integrating the above inequality on the interval $(0,x)$, we have
\begin{align}\label{xxx}
  e^{\sqrt{v_0}x}\varphi_-(x)-1
  \ge
  \frac{m_-+\sqrt{v_0}}{2\sqrt{v_0}}\left(e^{2\sqrt{v_0}x}-1\right)
  \ge \sqrt{\frac{v_0}{v_1}}\left(e^{2\sqrt{v_0}x}-1\right),
\end{align}
where the last inequality follows from $m_-\ge v_0/\sqrt{v_1}$ by \eqref{p3}.
From \eqref{xxx}, we get
\begin{align*}
  \varphi_-(x)\ge
  \sqrt{\frac{v_0}{v_1}}e^{\sqrt{v_0}x}
  +  \left( 1-\sqrt{\frac{v_0}{v_1}}\right)e^{-\sqrt{v_0}x}
  \ge  \sqrt{\frac{v_0}{v_1}}e^{\sqrt{v_0}x},
\end{align*}
which gives the lower bound estimate in \eqref{phi-} for $x>0$.

For the upper bound estimate in \eqref{phi-},
we define 
$y(x):=\varphi_-'(x)+\sqrt{v_1}\varphi_-(x)$ for $x\ge 0$.
Then, we have
\begin{align*}
  y'(x)=\varphi_-''(x)+\sqrt{v_1}\varphi_-'(x)=\sqrt{v_1}\varphi_-'(x)+V(x)\varphi_-(x)
  \le   \sqrt{v_1}\varphi_-'(x)+v_1\varphi_-(x)
  =\sqrt{v_1}y(x).
\end{align*}
Solving this differential inequality, we obtain
$y(x)\le y(0)e^{\sqrt{v_1}x}$. Since $y(0)=m_-+\sqrt{v_1}$,
this is equivalent to
\begin{align*}
  \varphi_-'(x)+\sqrt{v_1}\varphi_-(x)\le (m_-+\sqrt{v_1})e^{\sqrt{v_1}x}.
\end{align*}
We solve this differential inequality as follows:
\begin{align*}
  \left( e^{\sqrt{v_1}x}\varphi_-(x)\right)'
  =e^{\sqrt{v_1}x}\left( \varphi_-'(x)+\sqrt{v_1}\varphi_-(x)\right)
  \le (m_-+\sqrt{v_1})e^{2\sqrt{v_1}x}
  =\left( \frac{m_-+\sqrt{v_1}}{2\sqrt{v_1}}e^{2\sqrt{v_1}x}\right)'.
\end{align*}
Integrating the above inequality on the interval $(0,x)$, we have
\begin{align}\label{yyy}
  e^{\sqrt{v_1}x}\varphi_-(x)-1
  \le
  \frac{m_-+\sqrt{v_1}}{2\sqrt{v_1}}\left(e^{2\sqrt{v_1}x}-1\right)
\le  \sqrt{\frac{v_1}{v_0}}\left(e^{2\sqrt{v_1}x}-1\right),
\end{align}
where the last inequlity follows from $m_-\le v_1/\sqrt{v_0}$ by \eqref{p3}.
From \eqref{yyy}, we get
\begin{align*}
    \varphi_-(x)\le
  \sqrt{\frac{v_1}{v_0}}e^{\sqrt{v_1}x}
  + \left( 1-\sqrt{\frac{v_1}{v_0}}\right)e^{-\sqrt{v_1}x}
  \le   \sqrt{\frac{v_1}{v_0}}e^{\sqrt{v_1}x},
\end{align*}
which gives the upper bound estimate in \eqref{phi-} for $x>0$.
The estimate \eqref{phi+} can be proven smilarly,
and we omit it.

Lastly, we prove the uniqueness of the boundary value problems \eqref{ivp}.
We prove it for $\varphi_-$, since we can prove it for $\varphi_+$ similarly.
Let $\varphi_-^1$ and $\varphi_-^2$ be two solutions of \eqref{ivp}.
Setting $\varphi_0:=\varphi_-^1-\varphi_-^2$, it satisfies
\begin{align*}
  \begin{cases}
    \varphi_0''(x)=V(x)\varphi_0(x)\qquad(\text{a.e. }x\in\R),\\
    \varphi_0(0)=0,\\
    \DS{\lim_{x\to - \infty}\varphi_0(x)=0.}
  \end{cases}
\end{align*}
  Since a pair of $\varphi_+$ and $\varphi_-$ consists a basis
  of the solution space of the second order linear ODE $\varphi''+V(x)\varphi=0$,
  there exist $c_\pm\in\R$ such that
  $\varphi_0(x)=c_+\varphi_+(x)+c_-\varphi_-(x)$ holds for a.e. $x\in\R$.
  The condition $\lim_{x\to - \infty}\varphi_0(x)=0$ implies that $c_+=0$
  and $c_-=c_-\varphi_-(0)=\varphi_0(0)=0$ follows.
  Hence, the solution of \eqref{ivp} is unique.
\qed

\section{Second minimization step}\label{subsec:2.3}
\setcounter{equation}{0}

Based on \eqref{infF}, this section considers the second minimization problem
$\inf_{a\in\R}F(a)$.
The following theorem gives the derivative of the function $F(a)$.
\begin{Th}\label{FFF}
We suppose \eqref{bddV} and define $F(a)$ by \eqref{Fa}.
Then, $F\in W^{2,\infty}_{loc}(\R)$ and
the following estimates hold:
\begin{align}
&F(a)=u_a'(a_-)-u_a'(a_+)\quad (a\in\R),\label{F} \\
&F'(a)=|u_a'(a_+)|^2-|u_a'(a_-)|^2\quad (a\in\R),\label{F'}\\
  &F''(a)=2F(a)\left(|u_a'(a_+)|^2+u_a'(a_+)u_a'(a_-)+|u_a'(a_-)|^2-V(a)\right)
  \quad (\mbox{a.e. } a\in\R),\label{F''}
\end{align}
where $u_a'(a_+)$ and $u_a'(a_-)$ are the right and left derivatives of $u_a(x)$
at $x=a$, respectively.
Moreover, if $V$ is continuous on $\R$, $F\in C^2(\R)$ holds.
\end{Th}

\proof
From Theorem~\ref{phipm}, we have
\begin{align*}
F(a)
&=\int_{-\infty}^a\left(|u_a'(x)|^2+V(x)|u_a(x)|^2\right)dx
+\int_a^\infty\left(|u_a'(x)|^2+V(x)|u_a(x)|^2\right)dx\\
&=\int_{-\infty}^a\left(|u_a'(x)|^2+u_a''(x)u_a(x)\right)dx
+\int_a^\infty\left(|u_a'(x)|^2+u_a''(x)u_a(x)\right)dx\\
&=\int_{-\infty}^a\left( u_a(x)u_a'(x)\right)'dx
+\int_a^\infty \left(u_a(x)u_a'(x)\right)'dx\\
&=\Big[ u_a(x)u_a'(x)\Big]_{-\infty}^a
+\Big[ u_a(x)u_a'(x)\Big]_a^{\infty}\\
&=u_a(a)u_a'(a_-)-u_a(a)u_a'(a_+)=u_a'(a_-)-u_a'(a_+).
\end{align*}
From this formula and \eqref{ua'}, we also obtain $F\in W^{1,\infty}_{loc}(\R)$. 
Using Theorem~\ref{phipm}, we obtain
\begin{align*}
  \frac{d}{da}\left( u_a'(a_\pm)\right)
  =  \frac{d}{da}\left( \frac{\varphi_\pm'(a)}{\varphi_\pm (a)}\right)
  =  \frac{\varphi_\pm''(a)\varphi_\pm (a)-|\varphi_\pm'(a)|^2}{|\varphi_\pm (a)|^2}
  = V(a) -|u_a'(a_\pm)|^2,
\end{align*}
and \eqref{F'} follows.
We also have
\begin{align*}
  F'(a)=
  \left( u_a'(a_+)-u_a'(a_-)\right)  \left( u_a'(a_+)+u_a'(a_-)\right)
  =-F(a) \left( u_a'(a_+)+u_a'(a_-)\right).
\end{align*}
This formula implies $F\in W^{2,\infty}_{loc}(\R)$. 
Then, differentiating $F'(a)$, we obtain
\begin{align}
  F''(a)
  &=-F'(a) \left( u_a'(a_+)+u_a'(a_-)\right)-F(a)
  \left(2V(a)-|u_a'(a_+)|^2-|u_a'(a_-)|^2\right)\notag\\
&=2F(a)\left(|u_a'(a_+)|^2+u_a'(a_+)u_a'(a_-)+|u_a'(a_-)|^2-V(a)\right),\label{F2}
\end{align}
for a.e. $a\in \R$.
Moreover, if $V\in C^0(\R)$, $F\in C^2(\R)$ holds from \eqref{ua'} and \eqref{F2}.
\qed

\begin{Th}\label{mV-th}
Under the condition \eqref{bddV}, it holds that
\begin{align*}
  2\frac{v_0}{\sqrt{v_1}}\le m(V) \le 2\frac{v_1}{\sqrt{v_0}}.
\end{align*}
\end{Th}
\proof
From the estimate \eqref{p3}, we have
$v_0/\sqrt{v_1}\le \mp u_a'(a_\pm) \le v_1/\sqrt{v_0}$ for $a\in \R$.
Then, from \eqref{F}, 
$2v_0/\sqrt{v_1}\le F(a) \le 2v_1/\sqrt{v_0}$ holds.
Since $m(V)=\inf_{a\in\R}F(a)$, we conclude the assertion of the theorem.
\qed

When $V\in C^0(\R)$, if $F(x)$ has a local minimum at $x=a$, then
the following condition holds:
\begin{align}\label{lm}
  F'(a)=0,\qquad F''(a)\ge 0.
\end{align}
So, we define
$$N(V ):=\left\{a\in\R\;;\;F'(a)=0,F''(a)\ge0\right\}.$$
Then, as a consequence of Theorem~\ref{mV-th}, we obtain 
the following theorem.
\begin{Th}\label{ppp}
We suppose \eqref{bddV} and $V\in C^0(\R)$.   
Then,
\begin{align}\label{mVcond}
m(V)=\min\left(\liminf\limits_{|a|\to\infty}F(a),
\inf\limits_{a\in N(V)}F(a)\right),
\end{align}
\begin{align}\label{MVcond}
  M(V)=\left\{cu_a\;;\;c\in\R\setm\{0\},~F(a)=m(V)\right\}.
  \end{align}
In particular, $M(V )=\emptyset$ holds if and only if
\begin{align}\label{MVempty}
\liminf\limits_{|a|\to\infty}F(a)<F(b)\quad\mbox{for all}\quad b\in N(V).
\end{align}
\end{Th}
\proof
Since $V\in C^0(\R)$, $F\in C^2(\R)$ follows from Theorem~\ref{FFF}.
Then, \eqref{mVcond} holds from Proposition~\ref{decomposition}.
We set
\[
\tilde{M}(V):=\left\{cu_a\;;\;c\in\R\setm\{0\},~F(a)=m(V)\right\}.
\]
If $M(V)\neq \emptyset$, for any $u\in M(V)$,
we choose $a\in\R$ such that $|u(a)|=\| u\|_\infty$
and define $\tilde{u}(x):=u(x)/u(a)$.
Since $\tilde{u}\in K_a$, we obtain
\[
m(V)\le F(a)\le I(\tilde{u};V)=R(\tilde{u};V)=R(u;V)=m(V).
\]
Hence, $I(\tilde{u};V)=F(a)=m(V)$ follows.
From Theorem~\ref{th1}, we have $\tilde{u}=u_a$ and it implies $M(V)=\tilde{M}(V)$.
It is also obvious that $\tilde{M}(V)=\emptyset$ if $M(V)=\emptyset$.
So, we obtain \eqref{MVcond}.

When $M(V)=\emptyset$, since $m(V)=\liminf\limits_{|a|\to\infty}F(a)$ 
and $m(V)<F(b)$ for all $b\in N(V)$, \eqref{MVempty} holds.
  Conversely, \eqref{MVempty} also implies $M(V)=\emptyset$, since
    $m(V)\le \liminf\limits_{|a|\to\infty}F(a)$.
\qed

We immediately obtain characterizations of $m(V)$ and $M(V)$ for
the case of constant potential as a corollary of Theorem~\ref{ppp}.
\begin{Cor}\label{constantV}
Let $V $ be a positive constant:
$V >0$.
Then it holds that
\begin{align}\label{MV}
m(V)=2\sqrt{V},
~~~~~
M(V)=\{
cu_a;~c\in\R\setm\{0\},~a\in \R
\},
\end{align}
where $u_a(x)=e^{-\sqrt{V}|x-a|}$.
\end{Cor}
\proof
Since $v_0=v_1=V$, from Theorem~\ref{th1} and Theorem~\ref{mV-th},
we obtain $u_a(x)=e^{-\sqrt{V}|x-a|}$ and $m(V)=2\sqrt{V}$, respectively.
As $F(a)=\| u_a\|_V^2$ is constant, $N(V)=\R$ and \eqref{MV} follow from
Theorem~\ref{ppp}.
\qed

We also give characterizations of $m(V)$ and $M(V)$ for
the case of nondecreasing potential as follows.
\begin{Th}\label{nondecreasingV}
We suppose \eqref{bddV} and $V$ is a non-decreasing function with
\[
\lim_{x\to -\infty}V(x)=v_0>0,\quad \lim_{x\to \infty}V(x)=v_1.
\]
Then, it holds that
\[
m(V)=\lim_{a\to -\infty}F(a) =2\sqrt{v_0}.
\]
Furthermore, if $v_0<v_1$, then $F$ is a strictly increasing function
and $M(V)=\emptyset$.
\end{Th}
\proof
For $u\in H^1(\R)\setminus \{0\}$, since $V(x)\ge v_0$, we have
$R(u;V)\ge R(u;v_0)\ge m(v_0)$.
Taking the infimum of $R(u;V)$ with respect to $u$, we obtain
that $m(V)\ge m(v_0)$.

We put $v_a(x):=e^{-\sqrt{v_0}|x-a|}$ for $a\in\R$.
We also obtain
\begin{align*}
 m(v_0)&\le m(V)\le  F(a)\le I(v_a;V)\\
  &=I(v_a;v_0)+\left(I(v_a;V)-I(v_a;v_0)\right)\\
  & =m(v_0)+\int_{-\infty}^\infty (V(x)-v_0)e^{-2\sqrt{v_0}|x-a|}dx.
\end{align*}
Using the following estimate for $R>0$ and $a\in\R$, 
\begin{align*}
  0&< \int_{-\infty}^\infty (V(x)-v_0)e^{-2\sqrt{v_0}|x-a|}dx\\
  &=\int_{-\infty}^{a+R} (V(x)-v_0)e^{-2\sqrt{v_0}|x-a|}dx
  +\int_{a+R}^\infty (V(x)-v_0)e^{-2\sqrt{v_0}|x-a|}dx\\
  &\le \sup_{x\le a+R}(V(x)-v_0)
  \int_{-\infty}^\infty e^{-2\sqrt{v_0}|x-a|}dx
  +(v_1-v_0)\int_{a+R}^\infty e^{-2\sqrt{v_0}|x-a|}dx\\
  &= \frac{1}{\sqrt{v_0}}\sup_{x\le a+R}(V(x)-v_0)
  +\frac{v_1-v_0}{2\sqrt{v_0}}e^{-2\sqrt{v_0}R},
\end{align*}
we obtain
\[
\lim_{a\to -\infty}\int_{-\infty}^\infty (V(x)-v_0)e^{-2\sqrt{v_0}|x-a|}dx=0.
\]
and $m(v_0)\le m(V)\le \lim_{a\to -\infty}F(a) = m(v_0)$.
Hence, 
it holds that
\[
m(V)=\lim_{a\to -\infty}F(a) = m(v_0)=2\sqrt{v_0}.
\]

Furthermore, for $b<a$, we set $u_a\in K_a$ by \eqref{uaua}
and define $w_b\in K_b$ by $w_b(x):=u_a(x-b+a)$.
Then, if $v_0<v_1$, we have
\begin{align*}
  F(b)-F(a)
 \le I(w_b;V)-I(u_a;V)
 =\int_{-\infty}^\infty (V(x+b-a)-V(x)) |u_a(x)|^2dx <0.
\end{align*}
This estimate implies that $M(V)=\emptyset$, if $v_0<v_1$.
\qed

\section{A representation formula of the potential and its application}\label{subsec:2.4}
\setcounter{equation}{0}

We study the condition \eqref{lm} to investigate $m(V)$ and $M(V)$ in detail.
Surprisingly enough, we are able to express the {\it given} potential $V$
in terms of $\varphi_\pm$,
where $\varphi_\pm$ are the functions defined by Theorem~\ref{phipm}.
In this section, we observe this fact and use it to construct
an example of $V$, for which $m(V)$ and $M(V)$ are obtained explicitly.
We set
\begin{align*}
&h_\pm(x):=\varphi_\pm(x)^2,\quad
  \ell_\pm(x):=\log \varphi_\pm(x),\\
  &H_+(x):=\int_{-\infty}^x\frac{1}{h_+(\xi)}d\xi,\quad
    H_-(x):=\int_x^\infty\frac{1}{h_-(\xi)}d\xi,
\end{align*}
where $H_\pm(x)$ are well defined from the estimate \eqref{phi+} and \eqref{phi-}.
Then, we have the following proposition. 
\begin{Prop}
  We suppose \eqref{bddV} and set $W:=\varphi_-'(0)-\varphi_+'(0)>0$.
  Then, the following equalities hold:
\begin{align}
  &u_x'(x_\pm)=\frac{\varphi_\pm'(x)}{\varphi_\pm(x)}=\ell_\pm'(x)\quad (x\in\R),\label{L1}\\
  &\varphi_\mp(x)=W\varphi_\pm (x)H_\pm(x)\quad (x\in\R),\label{L2}\\[5pt]
  &V(x)=\ell_\pm''(x)+\ell_\pm'(x)^2\quad (\mbox{a.e. } x\in\R),\label{VV}\\[5pt]
  &F(x)=\ell_-'(x)-\ell_+'(x)\quad (x\in\R),\label{L3}\\
  &F'(x)=-F(x)(\ell_+'(x)+\ell_-'(x))=-\frac{WF(x)}{\varphi_+(x)\varphi_-(x)}
  (h_\pm'(x)H_\pm(x)\pm 1)\quad (x\in \R),\label{L4}\\
  &F''(x)=-2F'(x)\ell_\mp'(x)
  -2F(x)\ell_\pm ''(x)
  \quad (\mbox{a.e. } x\in\R).\label{L5}
\end{align}
\end{Prop}
\proof
The equality \eqref{L1} immediately follows from \eqref{uax} and the definition of $\ell_\pm$.
Then, 
\eqref{L3} and the first equality of \eqref{L4} follow
from Theorem~\ref{FFF} and \eqref{L1}.
Also, from \eqref{L1}, we have
\begin{align*}
  \ell_\pm''(x)=
  \frac{\varphi_\pm''(x)\varphi_\pm(x)-\varphi_\pm'(x)^2}{\varphi_\pm(x)^2}
  = V(x)-\ell_\pm'(x)^2,
\end{align*}
for a.e. $x\in\R$,
which implies the formula \eqref{VV}.

We define the Wronskian, 
$W:=\varphi_+(x)\varphi_-'(x)-\varphi_+'(x)\varphi_-(x)$.
Then, since 
\[
W'(x):=\varphi_+(x)\varphi_-''(x)-\varphi_+''(x)\varphi_-(x)
=\varphi_+(x)(V(x)\varphi_-(x))-(V(x)\varphi_+(x))\varphi_-(x)=0,
\]
holds for a.e. $x\in \R$, we obtain that $W$ is a constant
and $W=\varphi_-'(0)-\varphi_+'(0)>0$ holds.
Then, we obtain
\begin{align*}
  \frac{d}{dx}\left( \frac{\varphi_\mp(x)}{\varphi_\pm(x)}\right)
  =\frac{\varphi_\mp'(x)\varphi_\pm(x)-\varphi_\mp(x)\varphi_\pm'(x)}{\varphi_\pm(x)^2}
  =\frac{\pm W}{h_\pm(x)}=WH_\pm'(x).
\end{align*}
Integrating this equality,
we obtain $\varphi_\mp(x)/\varphi_\pm(x)=WH_\pm(x)$,
since $\lim_{x\to\mp\infty}H_\pm(x)=0$ and
$\lim_{x\to\mp\infty}\varphi_\mp(x)/\varphi_\pm(x)=0$.
Hence, \eqref{L2} holds.

From \eqref{L2}, we obtain $\varphi_+(x)\varphi_-(x)=Wh_\pm(x)H_\pm(x)$
and
\begin{align}\label{d1}
  \frac{d}{dx}\big(\varphi_+(x)\varphi_-(x)\big)
  =W\big( h_\pm'(x)H_\pm(x)+h_\pm(x)H_\pm'(x)\big)
  =W\big( h_\pm'(x)H_\pm(x)\pm 1\big),
\end{align}
where we used the relation $H_\pm'(x)=\pm h_\pm(x)^{-1}$.
On the other hand, we have
\begin{align}
  \frac{d}{dx}\big(\varphi_+(x)\varphi_-(x)\big)
  &=\varphi_+'(x)\varphi_-(x)+\varphi_+(x)\varphi_-'(x)\notag\\
  &=\varphi_+(x)\varphi_-(x)\left(
  \frac{\varphi_+'(x)}{\varphi_+(x)}+\frac{\varphi_-'(x)}{\varphi_-(x)}
  \right)\notag\\
  &=\varphi_+(x)\varphi_-(x)\big(
  \ell_+'(x) + \ell_-'(x)\big).\label{d2}
\end{align}
From \eqref{d1} and \eqref{d2}, we obtain
\begin{align*}
 \ell_+'(x) + \ell_-'(x)= \frac{W(h_\pm'(x)H_\pm(x)\pm 1)}{\varphi_+(x)\varphi_-(x)},
\end{align*}
which implies the second equality of \eqref{L4}.

Lastly, we show \eqref{L5}.
Differentiating \eqref{L3}, we have
$F'(x)=\ell_-''(x)-\ell_+''(x)$ for a.e. $x\in\R$.
So, it is equivalent to
\begin{align}\label{ellell}
  \ell_+''(x)+\ell_-''(x)=2\ell_\pm''(x)\pm F'(x)
  \quad(\mbox{a.e.}\quad x\in\R).
\end{align}
Hence, from the first equality of \eqref{L4} and \eqref{ellell},
we obtain \eqref{L5} as follows:
\begin{align*}
  F''(x)&=-F'(x)(\ell_+'(x)+\ell_-'(x))-F(x)(\ell_+''(x)+\ell_-''(x))\\
  &=-F'(x)(\ell_+'(x)+\ell_-'(x))-F(x)(2\ell_\pm''(x)\pm F'(x))\\
  &=-F'(x)(\ell_+'(x)+\ell_-'(x)\pm F(x))-2F(x)\ell_\pm''(x)\\
  &=-2F'(x)\ell_\mp '(x)-2F(x)\ell_\pm''(x).
\end{align*}
\qed

Thus, the necessary and sufficient condition for satisfying
condition \eqref{lm} is obtained as the following theorem.
\begin{Th}\label{study4}
We suppose \eqref{bddV} and $V\in C^0(\R)$.   
Then, each of the following conditions is equivalent to \eqref{lm}:
\begin{enumerate}
 \item $-\ell_+'(a)=\ell_-'(a)\geq \sqrt{V (a)}.$
 \item $h_+'(a)H_+(a)=-1
        \quad\mbox{and}\quad\ell_+''(a)\leq 0.$
 \item $h_-'(a)H_-(a)=1
        \quad\mbox{and}\quad\ell_-''(a)\leq 0.$
\end{enumerate}
\end{Th}
\noindent
\proof
From \eqref{F''} and \eqref{L1}, 
we obtain
\begin{align}\label{F''2}
  F''(x)=2F(x)\ell_+'(x)(\ell_+'(x)+\ell_-'(x))+
  2F(x)(\ell_-'(x)^2-V(x)).
\end{align}
From \eqref{F''2} and the first equality of \eqref{L4},
we conclude that condition~1 is equivalent to \eqref{lm}.
The equivalency of condition~2 or 3 to \eqref{lm} is
shown from \eqref{L4} and \eqref{L5}.
\qed

Using the new criterion which we obtained in Theorem~\ref{study4},
we give a nontrivial example of the bounded inhomogeneous potential 
for which the Sobolev best constant and function are explicitly specified.

\begin{Ex}
{\rm
For $A>0$ and $B>0$, we set
\begin{align}\label{vp+}
\varphi_+(x):=\frac{Ae^{-Bx}}{\sqrt{x^2+A^2}}.
\end{align}
Then, $\varphi_+(x)>0$ and we have
\begin{align}
\ell_+(x):=&\log \varphi_+(x)=\log A -Bx -\frac{1}{2}\log (x^2+A^2),\notag\\
\ell_+'(x)=&-B-\frac{x}{x^2+A^2}=-\frac{B(x^2+A^2)+x}{x^2+A^2},\notag\\
\ell_+''(x)=&-\frac{(x^2+A^2)-x(2x)}{(x^2+A^2)^2}
=\frac{x^2-A^2}{(x^2+A^2)^2}.\label{l+''}
\end{align}
Then, from \eqref{VV}, we obtain
\[
V (x):=\ell_+''(x)+\ell_+'(x)^2=\frac{(B(x^2+A^2)+x)^2+x^2-A^2}{(x^2+A^2)^2},
\]
and $V\in L^\infty(\R)$ holds.
We choose $A$ and $B$ as $V$ satisfies \eqref{bddV}.
Since
\begin{align*}
  V (x)&=B^2+\frac{2Bx}{x^2+A^2}+\frac{2x^2-A^2}{(x^2+A^2)^2}
  \ge B^2+\min_{x\in\R}\frac{2Bx}{x^2+A^2}+\min_{x\in\R}\frac{2x^2-A^2}{(x^2+A^2)^2}\\
  &=B^2 -\frac{B}{A}-\frac{1}{A^2}
  =\frac{(AB)^2-AB-1}{A^2},
\end{align*}
if $AB>(1+\sqrt{5})/2$, then $V (x)\ge (A^2B^2-AB-1)A^{-2}>0$ holds.
It is easy to check that $\varphi_+$ defined by \eqref{vp+}
satisfies the condition \eqref{ivp} for the above $V(x)$.

Since
\begin{align*}
&h_+(x)=\frac{A^2e^{-2Bx}}{x^2+A^2},
\quad
h_+'(x)=-\frac{2A^2\{ B(x^2+A^2)+x\}e^{-2Bx}}{(x^2+A^2)^2},\\
&H_+(x)=\frac{1}{A^2}\int_{-\infty}^x(x^2+A^2)e^{2Bx}\,dx
=\frac{(2B^2x^2-2Bx+2A^2B^2+1)e^{2Bx}}{4A^2B^3},
\end{align*}
using the formula \eqref{L2}, we have
\begin{align*}
\varphi_-(x)
&=W\varphi_+(x)H_+(x)\\
&=
\frac{W(2B^2x^2-2Bx+2A^2B^2+1)e^{Bx}}{4AB^3\sqrt{x^2+A^2}}\\
&=\frac{A(2B^2x^2-2Bx+2A^2B^2+1)e^{Bx}}{(2A^2B^2+1)\sqrt{x^2+A^2}},
\end{align*}
where the last equality follows from the condition $\varphi_-(0)=1$.
Then, we have
\begin{align*}
  \ell_-(x)&=\log \frac{A}{(2A^2B^2+1)}+
  \log (2B^2x^2-2Bx+2A^2B^2+1) -\frac{1}{2}\log (x^2+A^2)+Bx,\\
  \ell_-'(x)&=
  \frac{4B^2x-2B}{2B^2x^2-2Bx+2A^2B^2+1} -\frac{x}{x^2+A^2}+B.
\end{align*}

From \eqref{L3},
we obtain 
\[
\lim_{|a|\to \infty}F(a)=\lim_{|a|\to \infty}\big(\ell_-'(a)-\ell_+'(a)\big)=2B,
\]
since $\ell_\pm'(x)\rightarrow \mp B$ as $|x|\rightarrow \infty$.

Let us check the condition~2 of Theorem~\ref{study4}.
From \eqref{l+''}, the condition $\ell_+''(a)\leq 0$ is equivalent to $|a|\leq A$.
Using the identity:
\begin{align*}
  h_+'(x)H_+(x)
  &=-\frac{\{B(x^2+A^2)+x\}(2B^2x^2-2Bx+2A^2B^2+1)}{2B^3(x^2+A^2)^2}
  =\frac{Bx^2-x-A^2B}{2B^3(x^2+A^2)^2}-1,
\end{align*}
we obtain that the other condition $h_+'(a)H_+(a)=-1$
is equivalent to $Ba^2-a-A^2B=0$, namely
\[
a=a_1:=\frac{1- \sqrt{1+4A^2B^2}}{2B}\in (-A,0)
\quad
\mbox{or}\quad
a=a_2:=\frac{1+\sqrt{1+4A^2B^2}}{2B}>A.
\]
Hence we have $N(V )=\{a_1\}$. Since
\[
F(a_1)=-2\ell_+'(a_1)=2\left(B+\frac{a_1}{(a_1)^2+A^2}\right)
=2B\left(1-\frac{1}{\sqrt{1+4A^2B^2}}\right)
<2B=\lim_{|a|\rightarrow \infty}F(a),
\]
we conclude that 
\[
m(V )=2B\left(1-\frac{1}{\sqrt{1+4A^2B^2}}\right),
~~~~~
M(V )=\{cu_{a_1};~c\in\R\setm\{0\}\}.
\]
}
\end{Ex}

\section{Conclusion}
\setcounter{equation}{0}
We developed a new variational approach 
to determine the best constant of the
Sobolev-type embedding inequality in one dimension with a bounded inhomogeneous potential term.
We adopted two-step minimization method.
In the first minimization step, the Green function of the differential operator
$-\frac{d^2}{dx^2}+V(x)$
was captured as the minimizer.
We studied the fine properties
of the minimizer and the fundamental solutions of the ODE $-\frac{d^2}{dx^2}u+V(x)u=0$.
It enabled us to derive new necessary and sufficient conditions
for the local minimality in the second minimization step (Theorem~\ref{study4}).
Furthermore, as applications of our estimates, 
we constructed some concrete examples of the inhomogeneous 
potential $V(x)$ for which the best constant and function of the 
Sobolev-type embedding is exactly identified.

Since our approach provides a new tool to study the fine properties
of the best constant and function of the
Sobolev-type embedding with an inhomogeneous potential term,
further extension of our strategies 
is expected in future work. 
For example, the case the potential $V(x)$ is a periodic function
is interested in connection with the periodic crystal lattice, e.g., Kronig-Penney potential
\cite{K-P31}.
Extention of the Bloch theorem \cite{Kit04} for general
periodic potentials in our framework is also a challenging and worthy topic.
On the other hand, extending our method to the higher dimensional setting seems more complex. Nevertheless, it is a significant and exciting issue that we should address. Our approach in this paper is the first step toward that end.

~\\\noindent
{\bf Acknowledgment}:
This work was partially supported by JSPS KAKENHI Grant Nos. 20KK0058,
20H01812, and 20K03675.

\end{document}